\documentclass[12pt,reqno]{amsart}
\usepackage{amsmath}
\usepackage{amsfonts}
\usepackage{amssymb}
\usepackage{amsthm}
\usepackage{amscd}
\usepackage{a4wide}
\usepackage{enumerate}
\usepackage{dsfont} 

\usepackage{graphics}
\usepackage{eucal}
\usepackage{mathrsfs}

\usepackage{color}
\usepackage[pdftex,bookmarks,colorlinks,breaklinks]{hyperref}  

\hypersetup{linkcolor=red,citecolor=blue,filecolor=dullmagenta,urlcolor=darkblue} 

\newtheorem {theorem}    {Theorem}[section]
\newtheorem {lemma}      [theorem]    {Lemma}

\newtheorem {proposition}[theorem]    {Proposition}

\theoremstyle{definition}
\newtheorem {definition} [theorem]    {Definition}

\newcounter{AbcT}

\numberwithin{equation}{section}

\newcommand {\K} {{\mathbb K}}

\newcommand {\R} {{\mathbb R}}

\newcommand{\IGNORE}[1]{}

\renewcommand{\limsup}{\varlimsup}

 \DeclareMathOperator{\SL}{SL}

\DeclareMathOperator{\GL}{GL}

\newcommand{\frk}[2][n]{\mathfrak{#2}}

\newcommand{\mcl}[2][n]{\mathcal{#2}}

\newcommand\rank{\operatorname{rank}}
\newcommand\degree{\operatorname{deg}}
\newcommand\diameter{\operatorname{diameter}}

\newcommand{\di}{\operatorname{DI}_{\varepsilon}}

\begin{document}
\title{Dirichlet's theorem in function fields}

\begin{abstract}
We study metric Diophantine approximation in local fields of positive characteristic. Specifically, we study the problem of improving Dirichlet's theorem in Diophantine approximation and prove very general results in this context. 
\end{abstract}
\thanks{Ghosh is supported by an ISF-UGC grant.}

\author{Arijit Ganguly}
\author{Anish Ghosh}

\address{School of Mathematics, Tata Institute of Fundamental Research, Mumbai, 400005, India}
\email{arimath@math.tifr.res.in, ghosh@math.tifr.res.in}
\maketitle
\tableofcontents

\section{Introduction}

\subsection{The set up}

\noindent Let $p$ be a prime and $q:= p^r$, where $r\in \mathbb{N}$, let $\mathbb{F}_q$ be the finite field of $q$ elements and consider the field of rational functions $\mathbb{F}_{q}(T)$. 
We define a function $|\cdot|: \mathbb{F}_{q}(T) \longrightarrow \mathbb{R}_{\geq 0}$ as follows. 
\[ |0|:= 0\,\,  \text{ and} \,\, \left|\frac{P}{Q}\right|:= e^{\displaystyle \degree P- \degree Q}
\text{ \,\,\,for all nonzero } P, Q\in \mathbb{F}_{q}[T]\,.\] 
Clearly $|\cdot|$ is nontrivial,  non-archimedian and a discrete absolute value  
in $\mathbb{F}_{q}(T)$. This absolute value gives rise to a metric in $\mathbb{F}_{q}(T)$. \\

The completion field of $\mathbb{F}_{q}(T)$ with respect to this valuation is $\mathbb{F}_{q}((T^{-1}))$, the field of Laurent series 
over $\mathbb{F}_{q}$. The absolute value on $\mathbb{F}_{q}((T^{-1}))$, which we again denote by $|\cdot |$, is given as follows. 
Let $a \in \mathbb{F}_{q}((T^{-1}))$. For $a=0$, 
define $|a|=0$. If $a \neq 0$, then we can write 
$$a=\displaystyle \sum_{k\leq k_{0}} a_k T^{k}\,\,\mbox{where}\,\,\,\,k_0 \in \mathbb{Z},\,a_k\in \mathbb{F}_{q}\,\,\mbox{and}\,\, a_{k_0}\neq 0\,. $$
\noindent We define $k_0$ as the \textit{degree} of $a$, which will be denoted by $\degree a$,  and     
$|a|:= e^{\degree a}$. This clearly extends the absolute 
value $|\cdot|$ of $\mathbb{F}_{q}(T)$ to $\mathbb{F}_{q}((T^{-1}))$ and moreover, 
the extension remains non-archimedian and discrete like earlier. Let $\Lambda$ and $F$ 
denote $\mathbb{F}_{q}[T]$ and $\mathbb{F}_{q}((T^{-1}))$ respectively from now on. It is obvious that 
$\Lambda$ is discrete in  $F$. For any $n\in \mathbb{N}$, $F^n$ is throughout assumed to be equipped 
with the \textit{supremum norm} which is defined as follows
\[||\mathbf{x}||:= \displaystyle \max_{1\leq i\leq n} |x_i|\text{ \,\,for all \,} \mathbf{x}=(x_1,x_2,...,x_n)\in F^{n}\,,\]

\noindent and with the topology induced by this norm. Clearly $\Lambda^n$ is discrete in $F^n$. Since the topology on $F^n$ considered here 
 is the usual product topology on $F^n$, it follows that  $F^n$ is locally compact as $F$ is locally compact. We shall also
 fix a Haar measure $\lambda$ on $F$.\\

\noindent In this paper, we study analogues of Dirichlet's theorem in Diophantine approximation and its improvability for vectors in  $F^n$. An 
analogue of Dirchlet's theorem for local fields of positive characteristic  can be formulated as in the following:

\begin{theorem}\label{thm:D}
Let $t$ be a nonnegative integer. For $\mathbf{y}:=(y_1,y_2,...,y_n)\in F^n$, there exist $q\in \Lambda \setminus \{0\}$ and $p\in \Lambda$ 
such that
\[|y_1q_1+ y_2q_2+\cdot\cdot\cdot+y_nq_n-p|\,\textless \, \frac{1}{e^{nt}}\text{ \,\,and } \displaystyle 
\max_{1\leq j\leq n} |q_j|\leq e^t\,.\]
\end{theorem}

\noindent The theorem above is clearly well known and Diophantine approximation in the context of local fields of 
positive characteristic has been extensively studied of late. We refer the reader to \cite{deM} for survey and \cite{AGP, Kris, Las1, Las2} for more 
recent results. Indeed, the geometry of numbers, which can be used to prove Dirichlet's theorem was developed in the context of function fields by Mahler \cite{Mahler} as early as the 1940's. However, since we could not find a specific proof of the above in the literature, and in the interest of readability, we provide a proof  in 
Section \ref{section:review}. In fact we have proved there a stronger, namely a multiplicative statement (see Theorem \ref{thm:1.1}). There are many interesting parallels and contrasts between the theory of Diophantine approximation over the real numbers and in positive characteristic. Many results hold in both settings, the main result of the paper being one such while there are some striking exceptions. For instance the theory of \emph{badly approximable} numbers and vectors in positive characteristic offers several surprises: there is no analogue of Roth's theorem, provided that the base field is finite, which we assume throughout this paper. We refer the reader to \cite{AB} for other results in this vein. 

\subsection{Improving Dirichlet's theorem}Following Kleinbock-Weiss \cite{KW1}, the notion of 
``Dirichlet improvability" can now be introduced as follows. Let $0\,\textless\,\varepsilon \leq \frac{1}{e}$. A vector $\mathbf{y}:=(y_1,y_2,...,y_n)\in F^n$ is said to be \textit{Dirichlet $\varepsilon$-improvable}
if there is some $t_0\,\textgreater \,0$ such that for any choice of $n$ and nonnegative 
integers $t_1,t_2,...,t_n$ with $\max\{t,t_1,...,t_n\}\,\textgreater\,t_0$, where $t=t_1+
t_2+\,\cdot\,\cdot\,\cdot\,+t_n$, one can always 
find nonzero $(p,q_1,q_2,...,q_n)\in \Lambda \times \Lambda^n$ satisfying 
\[\displaystyle  |y_1q_1+ y_2q_2+\cdot\cdot\cdot+y_nq_n-p|\,\textless \, \frac{\varepsilon}{e^t} \text{ and } |q_j|\,\textless\, \varepsilon e^{t_j}\text{ for } j=1,2,...,n\,.
 \]

\noindent  Let $\di(n)$ denote the set of Dirichlet improvable vectors in $F^n$ or in $\R^n$, the context will make the field clear. Some remarks:\\
\begin{enumerate}
\item In the definition above, we have invoked the more general, multiplicative analogue of Dirichlet's thereorem, for which we provide a proof in Theorem \ref{thm:1.1}. The results of this paper are valid in this, stronger, setting.\\ 
\item This notion can be considered in greater generality, for systems of linear forms, as was done by Kleinbock and Weiss. We refer the reader to Definition \ref{defn:DI imp}.\\
\item Dirichlet's theorem can be formulated for global fields, i.e. one could consider number fields or finite extensions of positive characteristic fields. However, it is in general an open problem to determine the optimal constant in Dirichlet's theorem in this setting, without which of course, the question of improvement does not arise. There are some cases where the constant can be determined. For example in \cite{GR}, the theory of metric Diophantine approximation for certain \emph{imaginary} quadratic extensions of function fields was developed. In these fields, an analogue of Dirichlet's theorem with the same constant, i.e. $1$ holds, and it is plausible that the results of the present paper will work in that setting as well.\\

\end{enumerate}

\noindent We review briefly the state of the art on the question of improving Dirichlet's theorem 
in the context of real numbers. Davenport and Schmidt \cite{DS1, DS2} showed that the 
Lebesgue measure of $\di(n)$ is zero for every $\varepsilon < 1$. Starting with work of Mahler,  
the question of \emph{Diophantine approximation on manifolds} has received considerable attention. In this subject, 
one asks if Diophantine properties which are typical with respect to Lebesgue measure are also typical with respect to 
the push forward of Lebesgue measure via smooth maps. The starting point to this theory was a 
conjecture due to Mahler which asked if almost every point on the curve
\begin{equation}\label{def:veronese}
(x, x^2, \dots, x^n)
\end{equation}
\noindent is not very well approximable by rationals. Such maps (or measures) are referred to as \emph{extremal}. This 
conjecture was resolved by V. G. Sprindzhuk who in turn stated two generalisations 
of Mahler's conjecture which involved a \emph{nondegenerate} collection of functions replacing the map above. We refer to the work 
of Kleinbock-Margulis \cite{KM} where Sprindzhuk's conjectures are resolved, for the definitions. In a subsequent 
striking work, Kleinbock, Lindenstrauss and Weiss \cite{KLW} extended the results of \cite{KM} to a much wider 
class of measures, the so-called \emph{friendly} measures. This class includes push-forwards of Lebesgue measure as well as 
many other self similar measures including the uniform measure on the middle-third Cantor set. As regards improving Dirichlet's 
theorem for manifolds, in \cite{DS2}, Davenport and Schmidt showed that for any $\varepsilon < 4^{-1/3}$ the set of $x \in \R$ 
for which $(x, x^2) \in \di(2)$ has zero Lebesgue measure. Further results in this vein were obtained by Baker and by 
Bugeaud in \cite{Ba1, Ba2, Bu}. In \cite{KW1}, Kleinbock and Weiss proved several results in this direction and in particular 
showed the existence of $\varepsilon > 0$ such that for continuous, good and nonplanar maps $\mathbf{f}$ and Radon, Federer measures $\nu$, $\mathbf{f}_{*}(\nu)(\di(n)) = 0$. We will define all these terms later in the paper. In particular, this generalises the work of Baker and Bugeaud. The result that is obtained in \cite{KW1} holds for $\varepsilon$ which are quite a bit smaller than $1$ and to prove an analogous result for every $\varepsilon < 1$ remains an outstanding open problem. In the case of curves, N. Shah has resolved this problem. See \cite{Shah1} and also \cite{Shah2, Shah3} for related results.\\

In another direction, Kleinbock and Tomanov \cite{KT, KT-pre} established $S$-arithmetic 
analogues of Sprindzhuk's conjectures. In positive characteristic, Sprindzhuk \cite{Spr1} established the 
analogues of Mahler's conjecture, namely the extremality of the curve (\ref{def:veronese}) over $F^n$ and also proved other interesting results, including a transference principle interpolating between simultaneous Diophantine approximation and systems of linear forms. The analogues of Sprindzhuk's conjectures in positive characteristic were established by the second named author in \cite{G-pos}. However, the question of improving Dirichlet's theorem in positive characteristic has been completely open as far as we are aware. In the present paper, we study the question of Dirichlet improvability of vectors, maps and measures in positive characteristic.\\

\noindent Here is a special case of our main result, Theorem \ref{thm:1.2}. 
\begin{theorem}\label{thm:spcase}
Let $f_1,f_2,...,f_n$ be  polynomials so that $1,f_1,f_2,...,f_n$ are linearly independent over $F$. Fix some open set $U$ of $F$ and consider the map 
$\mathbf{f}(x)=(f_1(x),f_2(x),...,f_n(x))$ defined for all $x\in U$. Then there exists $\varepsilon_0\,\textgreater \,0$ such that 
whenever $\varepsilon\,\textless\,\varepsilon_0$,
 $\mathbf{f}(x)$ is not Dirichlet $\varepsilon$-improvable for $\lambda$ almost all $x\in U$.
 \end{theorem}
 
 \noindent Theorem \ref{thm:1.2}, the main result of this paper is far more general and holds for good, non-planar maps and Radon, Federer measures. It may be regarded as a positive characteristic version of 
Theorem 1.5 of Kleinbock and Weiss \cite{KW1}.  Since the statement of the general form of the Theorem is fairly technical, we have chosen to postpone it to later in the paper. The constant $\varepsilon$ can be estimated so the proof is ``effective" in that sense. However, it is likely to be far from optimal. We compute 
$\varepsilon_0$ in the special case $n=2$ and $f_i(x)=x^i$ for $i=1,2$ (see Section \ref{section:example}) as an example. Our 
proof proceeds along the lines of \cite{KW1} and the main tool is a quantitative non divergence result for certain maps in the 
space of unimodular lattices, which can be identified with the non compact quotient $\SL(n+1, F)/\SL(n+1, \Lambda)$. 
 
\section{Review of the classical theory}\label{section:review}
In this section, we provide a proof of Dirichlet's theorem in positive characteristic for completeness, and to aid the reader. In what follows, for  $k\in \mathbb{N}$, $\mathbb{Z}_{+}^{k}$ denotes the set of all 
$k$ tuples $(t_1,t_2,...,t_k)$ where each $t_i$ is a nonnegative integer . We prove the following:  
\begin{theorem}\label{thm:1.1} Let $m,n \in \mathbb{N}$, $k=m+n$ 
and $$\mathfrak{a}^{+}:= \{\mathbf{t}:=(t_1,t_2,...,t_k)\in \mathbb{Z}_{+}^{k}\,:\,\displaystyle \sum_{i=1}^{m} t_{i} =\displaystyle \sum_{j=1}^{n}t_{m+j}\}\,.$$ 
Consider $m$ linear forms $Y_{1},Y_{2},...,Y_{m}$  over $F$ in $n$ variables. Then for any $\mathbf{t}\in \mathfrak{a}^+$, 
there exist solutions $\mathbf{q}=(q_1,q_2,...,q_n)\in \Lambda^{n} \setminus \{\mathbf{0}\}$ and $\mathbf{p}=(p_1,p_2,...,p_m)\in \Lambda^m$ of the 
following system of inequalities 
\begin{equation}\label{eqn:1.1}\left \{ \begin{array}{rcl} |Y_{i}\mathbf{q}-p_i|\textless e^{-t_{i}} &\mbox{for} &i=1,2,...,m \\ |q_j|\leq e^{t_{m+j}}&\mbox{for} &j=1,2,...,n\,. \end{array} \right.  
\end{equation}
\end{theorem}
To prove this theorem, we first introduce the `polynomial part' and `fractional part' of a Laurent series. 
For any Laurent series $$a = \cdot\cdot\cdot+\frac{a_2}{T^2}+ \frac{a_1}{T}+(a_0+a_1T+a_2T^2+\cdot\cdot\cdot+a_kT^k)$$ in $F$, 
where $k \in \mathbb{Z},\,a_i\in \mathbb{F}_{q}$ and $a_{k}\neq 0$,  let us define the \textit{polynomial part} 
of $a$ as $$ a_0+a_1T+a_2T^2+\cdot\cdot\cdot+a_k T^k $$ if $k\geq 0$,  otherwise it is defined to be $0$;  and the \textit{fractional part} 
of $a$, denoted by $\langle a\rangle$, is defined as $$\alpha- \mbox{polynomial\,\,part}= \frac{a_1}{T}+\frac{a_2}{T^2}+\cdot\cdot\cdot\,. $$ \\

Now, let $a:=\frac{a_1}{T}+\frac{a_2}{T^2}+\cdot\cdot\cdot \in F$ and $\alpha:=\alpha_{0}+\alpha_{1}T+\alpha_{2}T^2+\cdot\cdot\cdot+\alpha_{k}T^k \in \Gamma\setminus \{0\}$ 
with degree $\leq k$, where $k\geq 0$ is an integer. Let us observe that, 
for any $s\in \mathbb{N}$, the coefficient 
of $\frac{1}{T^s}$ in $\alpha a$ 
is $$a_s \alpha_0+\cdot\cdot\cdot+a_{s+k}\alpha_{k}\,.$$ It follows that, for any $m\in \mathbb{N}$, $|\langle \alpha a\rangle |\textless \frac{1}{e^m}$ if and only if the  
system   $A\textit{\textbf{x}}=\textbf{0}$ of linear equations over $\mathbb{F}_{q}$, where  the coefficient matrix 
$$A:=\begin{bmatrix}a_1& a_2 & .&.&. & a_{k+1}\\a_2&
           a_3 & .&.&. & a_{k+2}\\.&.&\,&\,&\,&.\\.&.&\,&\,&\,&.\\.&.&\,&\,&\,&.\\a_m&a_{m+1} &.&.&.&a_{m+k}
         \end{bmatrix}\,,$$
  has  $(\alpha_0,\alpha_1,...,\alpha_k)$ as a nontrivial solution. \\
         
         Continuing along the same line, let us now take two Laurent series 
         $a=\frac{a_1}{T}+\frac{a_2}{T^2}+\cdot\cdot\cdot$, $b=\frac{b_1}{T}+\frac{b_2}{T^2}+\cdot\cdot\cdot,$ 
         and two nonzero polynomials $\alpha:=\alpha_{0}+\alpha_{1}T+\alpha_{2}T^2+\cdot\cdot\cdot+\alpha_{k}T^k$, $\beta:=\beta_{0}+\beta_{1}T+\beta_{2}T^2+\cdot\cdot\cdot+\beta_{l}T^l$, 
         with degree $\leq k,l$ respectively, where $k,l=0,1,2,...$\,\,. For any $s\in \mathbb{N}$, the coefficient 
         of $\frac{1}{T^s}$ in $\alpha a+\beta b$ is easily seen to be $$a_s \alpha_0+\cdot\cdot\cdot+a_{s+k}\alpha_{k}+b_s \beta_0+\cdot\cdot\cdot+b_{s+l}\beta_{l}\,.$$ 
         Therefore, for any $m\in \mathbb{N}$, $|\langle \alpha a+\beta b\rangle|\textless \frac{1}{e^m}$ if and only if $(\alpha_0,\alpha_1,...,\alpha_k,\beta_0,\beta_1,..., \beta_l)$ is a nontrivial solution of the following system
$$\begin{bmatrix}\Huge A &\huge B\end{bmatrix}
         \begin{bmatrix}\textit{\textbf{x}}\\ \textit{\textbf{y}}
         \end{bmatrix}=\textbf{0}\,,$$ where 
         $$A:= \begin{bmatrix}a_1& a_2 & .&.&. & a_{k+1}\\a_2&
           a_3 & .&.&. & a_{k+2}\\.&.&\,&\,&\,&.\\.&.&\,&\,&\,&.\\.&.&\,&\,&\,&.\\a_m&a_{m+1} &.&.&.&a_{m+k}
         \end{bmatrix}\,\,\mbox{and}\,\, B:=\begin{bmatrix}b_1& b_2 & .&.&. & b_{l+1}\\b_2&
           b_3 & .&.&. & b_{l+2}\\.&.&\,&\,&\,&.\\.&.&\,&\,&\,&.\\.&.&\,&\,&\,&.\\b_m&b_{m+1} &.&.&.&b_{m+l}
         \end{bmatrix}\,.$$ It is obvious that we can generalize this observation for any such $n$ Laurent series and nonzero polynomials. \\
         
       \indent Now we are ready to start the proof of Theorem \ref{thm:1.1}. Each $Y_i$, being a linear form over $F$ in $n$ variables, 
       must be of the form $$y_{i1}x_1+y_{i2}x_2+\cdot \cdot \cdot +y_{in}x_n\,,$$ for some $y_{ij}\in \mathbb{F}_q$ , $j=1,2,...,n$. It 
       suffices to consider the case $|y_{ij}|\,\textless\, 1$, i.e. the polynomial part of $y_{ij}$ is zero, for all $i=1,2,...,m$ 
       and $j=1,2,...,n$.  \\
       
       From the observations we made earlier, we see that each $y_{ij}$ gives rise to a matrix $M_{ij}$ having $t_i$ rows and $t_{m+j}+1$ columns 
and more importantly, the existence of solution of the system (\ref{eqn:1.1}) is equivalent to the existence of nontrivial solutions
of the following system of linear equations over $\mathbb{F}_q$
 $$\begin{bmatrix}M_{11}& M_{12} & .&.&. & M_{1n}\\M_{21}&
           M_{22} & .&.&. & M_{2n}\\.&.&\,&\,&\,&.\\.&.&\,&\,&\,&.\\.&.&\,&\,&\,&.\\M_{m1}&M_{m2} &.&.&.&M_{mn}
         \end{bmatrix}
         \begin{bmatrix}\textit{\textbf{x}}_1\\ \textit{\textbf{x}}_2\\ .\\.\\.\\ \textit{\textbf{x}}_n
         \end{bmatrix}=\textbf{0}\,.
         $$

Clearly the above coefficient matrix has  $\displaystyle \sum_{i=1}^m t_i$  rows and $\displaystyle\sum_{j=1}^n (t_{m+j}+1)$ columns. 
As $ \sum_{i=1}^{m} t_{i} = \sum_{j=1}^{n}t_{m+j}$, we see that 
the matrix has more columns than rows and hence nontrivial solution exists. This completes our proof. $\Box$\\
\section{The main theorem}
\indent We shall now introduce the notion of ``Dirichlet improvability'' in a greater generality. Let $\mathfrak{a}^{+}$ be as given in Theorem
\ref{thm:1.1},  $\mathcal{T}$ be an unbounded subset
of  $\mathfrak{a}^{+}$ and  $0\,\textless \,\varepsilon\leq \frac{1}{e}$. 
\begin{definition}\label{defn:DI imp} For a system of linear forms  $Y_{1},Y_{2},...,Y_{m}$  over $F$ in $n$ variables, 
we say that \textit{DT can 
be $\varepsilon$-improved along $\mathcal{T}$}, or we use the notation $Y\in DI_{\varepsilon}(\mathcal{T})$, where $Y$ is 
the $m\times n$ matrix having
$Y_i$ as the $i$ th row for each $i$,  if there exists $t_0\,\textgreater \,0$ such that for every 
$\mathbf{t}:=(t_1,t_2,...,t_k)\in \mathcal{T}$ with 
$||\mathbf{t}||\,\textgreater\, t_0$ the following system admits nontrivial  solutions $(\mathbf{p},\mathbf{q})\in \Lambda^m \times \Lambda^n$ : 
        \begin{equation}\label{eqn:1.4} \left \{\begin{array}{rcl} |Y_{i}\mathbf{q}-p_i|\,\textless \,\frac{\varepsilon}{e^{t_i}} &\mbox{for} &i=1,2,...,m 
        \\ |q_j|\,\textless\,\varepsilon e^{t_{m+j}}&\mbox{for} &j=1,2,...,n\,. \end{array} \right.\end{equation}   
 In particular, a vector $\mathbf{y}:=(y_1,y_2,...,y_n)\in F^n$ is said to be \textit{Dirichlet $\varepsilon$-improvable along $\mathcal{T}$} 
if the corresponding row matrix $[y_1\,\,y_2\,\,\cdot\,\,\cdot\,\,\cdot\,\,y_n]\in DI_{\varepsilon}(\mathcal{T})$.
\end{definition}

Exactly similar to that shown in \cite{KW1},  here also we want to prove that if  $\varepsilon\,\textgreater \,0$ is sufficiently small 
and an unbounded subset $\mathcal{T}$ of $\mathfrak{a}^{+}$ is chosen, the set of all  Dirichlet $\varepsilon$-improvable 
vectors along $\mathcal{T}$ is negligible. The setup here is \emph{multiplicative}, i.e. one studies Diophantine inequalities where Euclidean or supremum norm is replaced with the product of coordinates. The changed ``norm" introduces several complications and the subject of multiplicative Diophantine approximation is generally considered more difficult than its euclidean counterpart.\\

Before proceeding to our main theorem, we will recall the some terminology  
introduced in the papers of Kleinbock and Margulis, and Kleinbock, Lindenstrauss and Weiss and used in several subsequent works by many authors. The following is taken from  \S 1 and 2 of \cite{KT}. \\

For the sake of generality, we assume $X$ is a Besicovitch metric space, $U\subseteq X$ is  open, $\nu$ is a radon measure
on $X$, $(\mcl{F},|\cdot|)$ is a valued 
field and $f: X\longrightarrow\mcl{F}$ is a given function such that $|f|$ is measurable. For any $B\subseteq X$, we set
$$ ||f||_{\nu,B} := \displaystyle \sup_{x\in B\cap \text{ supp }(\nu)} |f(x)|.$$

\begin{definition}\label{defn:C,alpha}
 For $C,\alpha \textgreater\,0$, $f$ is said to be $(C,\alpha)-good$ on $U$ with respect to $\nu$ if for every ball $B\subseteq U$
 with center in $\text{supp }(\nu)$, one has
 \[\nu(\{x\in B: |f(x)|\,\textless \varepsilon\})\leq C\left(\frac{\varepsilon}{||f||_{\nu,B}}\right)^{\alpha} \nu(B)\,.\]
\end{definition}
The following properties are immediate from Definition \ref{defn:C,alpha}. 
\begin{lemma} \label{lem:C,alpha} Let $X,U,\nu, \mcl{F}, f, C,\alpha,$ be as given above. Then one has 
\begin{enumerate}[(i)]
   \item $f$ \text{ is } $(C,\alpha)-good$ \text{ on  }$U$ with respect to $\nu \Longleftrightarrow \text{ so is } |f|$.
   \item $f$ is $(C,\alpha)-good$ on $U$  with respect to $\nu$ $\Longrightarrow$ so is $c f$ for all $c \in \mcl{F}$.
   \item \label{item:sup} $\forall i\in I, f_i$ are $(C,\alpha)-good$ on $U$ with respect to $\nu$ and $\sup_{i\in I} |f_i|$ is measurable $\Longrightarrow$ so is $\sup_{i\in I} |f_i|$.
   \item $f$ is $(C,\alpha)-good$ on $U$ with respect to $\nu$ and $g :V \longrightarrow \mathbb{R}$ is a continuous function such 
   that $c_1\leq |\frac{f}{g}|\leq c_2$ 
 for some $c_1,c_2 \,\textgreater \,0\Longrightarrow g$ is $(C(\frac{c_2}{c_1})^{\alpha},\alpha)$ good on $U$ with respect to $\nu$.
 \item Let $C_2 \,\textgreater \,1$ and 
 $\alpha _2\,\textgreater\,0$. $f$ is $(C_1,\alpha_1)-good$ on $U$ with respect to $\nu$ and $C_1 \leq C_2, \alpha_2 \leq \alpha_1 \Longrightarrow f$ is
 $(C_2,\alpha_2)-good$ on $V$ with respect to $\nu$.
  \end{enumerate}
\end{lemma}  
 We say a  map $\mathbf{f}=(f_1,f_2,...,f_n)$ from $U$ to $\mcl{F}^n$, where $n\in \mathbb{N}$,
is $(C,\alpha)-good$ on $U$ with respect to $\nu$, or simply $(\mathbf{f},\nu)$ is 
$(C,\alpha)-good$ on $U$,  if every $\mcl{F}$-linear combination of $1,f_1,...,f_n$ 
is $(C,\alpha)-good$ on $U$ with respect to $\nu$. 

\begin{definition} Let $\mathbf{f}=(f_1,f_2,...,f_n)$ be a map from $U$ to $\mcl{F}^n$, where $n\in \mathbb{N}$. 
We say that $(\mathbf{f},\nu)$ is \emph{nonplanar} if for any ball $B\subseteq U$ with center in $\text{supp }(\nu)$, 
the restrictions of the functions  $1,f_1,...,f_n$ 
on $B\cap \text{ supp }(\nu)$ are linearly independent. 
 \end{definition}
In other words, $\mathbf{f}(B\cap \text{ supp }(\nu))$ is not contained in any affine subspace of $\mcl{F}^n$ 
for any ball $B\subseteq U$ with center in $\text{supp }(\nu)$.\\

For $m\in \mathbb{N}$ and a ball $B=B(x;r)\subseteq X$, where $x\in X$ and $r\,\textgreater\,0$, we shall use the notation
$3^mB$ to denote the ball $B(x;3^mr)$. 
\begin{definition}
 Let $D\,\textgreater\,0$. The measure $\nu$ is said to be $D-Federer$ on $U$ if for every ball $B$ with center 
 in $\text{supp }(\nu)$ such that $3B\subseteq U$, one has \[\frac{\nu(3B)}{\nu(B)}\leq D\,.\]
\end{definition}

\noindent We are now ready to state our main Theorem, which addresses improvements of Dirichlet's theorem in the multiplicative setting for good, nonplanar maps and Federer measures over local fields of positive characteristic.

\begin{theorem}\label{thm:1.2}
 For any $d,n\in \mathbb{N}$ and $C,\alpha,D\,\textgreater\, 0$ there exists $\varepsilon_0=\varepsilon_0(n,C,\alpha,D)$ satisfying
 the following: 
 whenever a radon  
 measure $\nu$
 on $F^d$, an open set $U$ of $F^d$ with $\nu(U)\,\textgreater \,0$ and $\nu$ is D-Federer on $U$, and  
 a continuous map $\mathbf{f}:U\longrightarrow F^n$ such that 
 $(\mathbf{f},\nu)$ is $(C,\alpha)-good$ and nonplanar is given then for any $\varepsilon < \varepsilon_{0}$, 
 
 \begin{equation}
  \mathbf{f}_{*}\nu(DI_{\varepsilon}(\mathcal{T}))=0 \,\,\text{ for any unbounded \,}\mathcal{T}\subseteq \mathfrak{a}^{+}\,.
 \end{equation}
\end{theorem}
We shall use the so called  ``quantitative nondivergence'', a  generalization of non-divergence of 
unipotent flows on homogeneous spaces, to prove our main theorem. Similarly to the approach adopted in \cite{KW1}, we will
first translate the property of a system of linear forms  over $F$ being Dirichlet improvable into certain recurrent properties of flows on 
some homogeneous space in the following section. \\
\section{The correspondence}\label{section:2}
Let $G:= \SL(k,F)$, $\Gamma := \SL(k, \Lambda)$ and $\pi$ be the quotient map $G\longrightarrow G/\Gamma$. $G$ acts on $G/\Gamma$  
by left translations via the rule $g\pi(h)=\pi(gh)$ for $g,h\in G$. For $Y\in M_{m\times n}(F)$, define $$\tau (Y):= \begin{bmatrix}
                                                                                         I_m& Y\\0& I_n
                                                                                        \end{bmatrix}\,\,\mbox{and}\,\,\overline{\tau}:=\pi \circ \tau\,,$$
                                                                                        where $I_l$ stands for the $l\times l$ identity matrix, $l\in \mathbb{N}$. 
                                                                                        Since $\Gamma$ is the stabilizer of $\Lambda^k$ under the transitive 
                                                                                        action of 
$G$ on the set of unimodular lattices in $F^k$, which is denoted by $\mathcal{L}_k(F)$, we can identify $G/\Gamma\simeq \mathcal{L}_k(F)$.
Thus $\overline{\tau}(Y)$ becomes identified with $$\{(Y\mathbf{q}-\mathbf{p},\mathbf{q})\,:\mathbf{p}\in \Lambda^m\,,\mathbf{q}\in \Lambda^n\}\,.$$

Now for $\varepsilon\textgreater 0$, let $K_{\varepsilon}$ denote the collection of all unimodular 
lattices in $F^k$ which contain nononzero vector of norm smaller than $\varepsilon$, that is, 
\begin{equation} \label{eqn:cpt} K_{\varepsilon}:=\pi (\{g\in G\,:\,||g\mathbf{v}||\geq \varepsilon \,\,\forall 
 \,\mathbf{v}\in \Lambda^k\setminus \{\mathbf{0}\}\})\,.
\end{equation}

Next, for $\mathbf{t}:=(t_1,t_2,...,t_k)\in \mathfrak{a}^{+}$, we associate the diagonal matrix 
$$g_{\mathbf{t}}:=\,\,\mbox{diag}\,(T^{t_1},.\,.\,.\,,T^{t_m},T^{-t_{m+1}},.\,.\,.\,,T^{-t_{k}})\in G\,.$$
Let us come to the relevance of defining the above objects. An immediate observation shows that, 
for given $\mathbf{t}\in \mathfrak{a}^{+}$, the system (\ref{eqn:1.4}) has nonzero polynomial solutions if and only if  
$$g_{\mathbf{t}}\overline{\tau}(Y)\notin K_{\varepsilon}\,.$$ Thus we have
\begin{proposition}\label{prop: 2.1}
 Let $0< \varepsilon \leq \frac{1}{e}$ and unbounded $\mathcal{T}\subseteq \mathfrak{a}^+$ be  given. Then for any $Y\in M_{m\times n}(F)$,
 $$Y\in DI_{\varepsilon}(\mathcal{T})\Longleftrightarrow g_{\mathbf{t}}\overline{\tau}(Y)\notin 
 K_{\varepsilon} \,\,\forall \mathbf{t}\in 
 \mathcal{T}\,\,\mbox{with} \,\,||\mathbf{t}||\gg 1\,,$$ or equivalently one has,
 $$
  DI_{\varepsilon}(\mathcal{T})= \displaystyle \bigcup_{n=1}^{\infty} \,\,\bigcap_{\mathbf{t}\in \mathcal{T}\,,\,||\mathbf{t}||\textgreater n}
  \{Y\in M_{m\times n}(F):\,g_{\mathbf{t}}\overline{\tau}(Y)\notin K_{\varepsilon}\}\,.
 $$

\end{proposition}
Hence, in view of the above proposition, it is clear that if in addition  a  radon
 measure $\nu$
 on $F^d$, an open set $U$ of $F^d$  and  a  map $F:U\longrightarrow M_{m\times n}(F)$ are given then to prove 
 $F_{*}\nu(DI_{\varepsilon}(\mathcal{T}))=\nu (F^{-1}(DI_{\varepsilon}(\mathcal{T}))=0$, it is enough to show 
  \begin{equation}\label{eqn:2.1} \nu \left(F^{-1} \left( \displaystyle \bigcap_{\mathbf{t}\in \mathcal{T}\,,\,||\mathbf{t}||\textgreater n}
  \{Y\in M_{m\times n}(F):\,g_{\mathbf{t}}\overline{\tau}(Y)\notin K_{\varepsilon}\}\right)\right) = 0 \end{equation} for all $n\in \mathbb{N}$. Suppose 
  now that 
  we have some $c \in (0,1)$ with the property that for any ball $B\subseteq U$ centered in $\mbox{supp}\,(\nu)$, there 
  exists $s\,\textgreater\, 0$ such that
  \begin{equation} \label{eqn:2.2}\displaystyle 
 \nu( B \cap F^{-1}(\{Y\,:\,g_{\mathbf{t}}\overline{\tau}(Y)\notin K_{\varepsilon}\}))= \nu 
 (\{\textbf{x}\in B\,:\,g_{\mathbf{t}}\overline{\tau}(F(\textbf{x}))\notin K_{\varepsilon}\})\leq c\nu (B) 
   \end{equation} holds for  any $\mathbf{t}\in \mathfrak{a}^+$ with $||\mathbf{t}||\geq s$.
 Then it is easy to see that, for any $n\in \mathbb{N}$ and 
  any ball $B\subseteq U$  centered in $\mbox{supp}\,(\nu)$,
\begin{equation} \label{eqn:2.3} \begin{array}{rcl}\displaystyle \frac{\nu \left(B \cap F^{-1} \left ( \displaystyle \bigcap_{\mathbf{t}\in \mathcal{T}\,,\,||\mathbf{t}||\textgreater n}
  \{Y\in M_{m\times n}(F):\,g_{\mathbf{t}}\overline{\tau}(Y)\notin K_{\varepsilon}\}\right)\right)}{\nu(B)}\\=
  \displaystyle \frac{\nu \left(\displaystyle \bigcap_{\mathbf{t}\in \mathcal{T}\,,\,||\mathbf{t}||\textgreater n}
  B\cap F^{-1}(\{Y\in M_{m\times n}(F):\,g_{\mathbf{t}}\overline{\tau}(Y)\notin K_{\varepsilon}\})\right)}{\nu(B)}& \leq 
  \displaystyle \frac{c\,\nu(B)}{\nu(B)}=c\textless 1\,.\end{array}\end{equation}

\noindent It follows that, for any given $n\in \mathbb{N}$, no $\mathbf{x} \in U\cap\mbox{supp}\,(\nu)$  is a point of 
density  of
the set  
$$F^{-1} \left(\displaystyle \bigcap_{\mathbf{t}\in \mathcal{T}\,,\,||\mathbf{t}||\textgreater n} 
\{Y\in M_{m\times n}(F):\,g_{\mathbf{t}}\overline{\tau}(Y)\notin K_{\varepsilon}\}\right)\,,$$
as (\ref{eqn:2.3}) holds true for any ball $B$ with $\mathbf{x}\in B\subseteq U$. Thus (\ref{eqn:2.1}) will be achieved in view of 
Theorem \ref{thm:3.1}.

\section{The proof of Theorem \ref{thm:1.2}}
As $F^d$ is locally compact, hausdorff and second countable, every open set is the union of some countable collection of 
compact subsets. Hence to prove the Theorem \ref{thm:1.2}, once correct $\varepsilon_0=\varepsilon_0(n,C,\alpha,D)$ is found,  
it suffices to show that for  all $\mathbf{y}\in U\cap\text{ supp }(\nu)$, 
there exists a ball $\mathfrak{B}\subseteq U$ containing $\mathbf{y}$ such that 
\begin{equation}\label{eqn:4.1} \nu(\mathfrak{B}\cap \mathbf{f}^{-1}(DI_{\varepsilon}(\mathcal{T})))=
\nu(\{\mathbf{x}\in \mathfrak{B}\,:\, \mathbf{f}(\mathbf{x})\in DI_{\varepsilon}(\mathcal{T})\})=0 \end{equation} for 
all $\varepsilon\,\textless\, \varepsilon_0$. 
From our discussion of Section \ref{section:2}, we see that (\ref{eqn:4.1}) is guaranteed as soon as we can show the 
existence of some $c\in (0,1)$ which satisfies the following: whenever a ball $B$ with center in $\mbox{supp}\,(\nu)$ is contained in $\mathfrak{B}$ then, 
\begin{equation}\label{eqn:suff}
 \mbox{there exists}\,\, s\,\textgreater \,0 \,\,\mbox{such that for all}\,\,\mathbf{t}\in \mathfrak{a}^+\,\
 \mbox{with}\,\,||\mathbf{t}||\geq s, (\ref{eqn:2.2})\,\,\mbox{holds}\,.
\end{equation}
Now the following proposition shows our way. 
\begin{proposition}\label{main prop} 
 For any $d,n\in \mathbb{N}$ and any $C,\alpha, D\,\textgreater \,0$ there exists $\tilde C=\tilde C(n,C, D)$ with the following property:\\
 
 \noindent whenever a ball $ B$ centered in $\mbox{supp}\,(\nu)$, a  radon measure  $\nu$ on $F^d$ which is D-Federer 
 on $\tilde B:= 3^{n+1} B$ 
 and a continuous map $\mathbf{f}:\tilde {B} \longrightarrow F^n$ are given so that 
 \begin{enumerate}[(i)]
  \item any $F$-linear combination of $1,f_1,\,.\,.\,.\,,f_n$ is $(C,\alpha)-good$ on $\tilde B$ with respect to $\nu$ and\,, 
 \item the restrictions of $1,f_1,\,.\,.\,.\,,f_n$ to $B\cap \mbox{supp}\,(\nu)$  are linearly independent 
 over $F$; \end{enumerate}
\noindent then we can find some $s\,\textgreater\, 0$ such that for all $\mathbf{t}\in \mathfrak{a}^+$ with $||\mathbf{t}||\geq s$ and 
any $\varepsilon\leq \frac{1}{e}$, one has 
\begin{equation}\label{eqn:prop}
\nu(\{\mathbf{x}\in B\,:\,g_{\mathbf{t}}\overline{\tau}(\mathbf{f}(\mathbf{x}))\notin K_{\varepsilon}\})\leq \tilde C \varepsilon ^{\alpha} \nu(B)\,.
\end{equation}
\end{proposition}
Theorem \ref{thm:1.2} follows easily from the Proposition \ref{main prop}. In fact, we first choose $ 0\,\textless\, \varepsilon_0\leq \frac{1}{e}$ 
so that $\displaystyle \tilde C \varepsilon_{0}  ^{\alpha}\,\textless \,1$. Clearly this $\varepsilon_0$ depends only on $(n, C,\alpha, D)$. 
Let $\mathbf{y}\in U\cap\mbox{supp}\,(\nu)$. Choose a ball $\mathfrak{B}$ such that 
$\mathbf{y}\in \mathfrak{B} \subseteq \tilde {\mathfrak{B}}:= 3^{n+1} \mathfrak{B} \subseteq U$. 
Now pick any 
ball $B\subseteq \mathfrak{B}$ having center in $\mbox{supp}\,(\nu)$ and consider the corresponding $\tilde B$.
Since $(\mathbf{f},\nu)$ is $(C,\alpha)-good$ and nonplanar, the conditions (i) and (ii) of Proposition \ref{main prop} hold here
immediately. Hence, 
if we set $c= \tilde C \varepsilon_{0}  ^{\alpha}$, the assertion (\ref{eqn:suff}) 
is immediate from Proposition \ref{main prop} whenever $0\,\textless\, \varepsilon \,\textless\, \varepsilon_0$. 
Thus the proof of Theorem \ref{thm:1.2} is complete. \,\,$\Box$\\

We now need to prove Proposition \ref{main prop}. We shall show this as a consequence of a more general 
result, namely the `Quantitative nondivergence  theorem'. All these will be discussed in Section \ref{section:5}.  

\section{Quantitative nondivergence and the proof of \\Proposition \ref{main prop}} \label{section:5}
We shall first recall the `Quantitative nondivergence theorem' in the most generality, as it is developed in \S 6 of \cite{KT}.
Finally,  we shall prove Proposition \ref{main prop} from this.

\subsection{Quantitative nondivergence}  
We start this subsection by assuming that $\mathcal{D}$ is an integral domain, $K$ is the field of quotients of $\mathcal{D}$ and 
$\mathcal{R}$ is a commutative ring containing $K$ as a subring. \\

Let $m\in \mathbb{N}$. If $\Delta$ is a $\mathcal{D}$-submodule of $\mathcal{R}^m$, let us 
denote by $K\Delta$ (respectively $\mathcal{R}\Delta$) its $K$- (respectively $\mathcal{R}$) linear span
inside $\mathcal{R}^m$. We use the notation $\rank(\Delta)$ to denote the rank of $\Delta$ which is defined as 
\[ \rank(\Delta):= \displaystyle \dim_{K} (K\Delta)\,.\]
For example $\rank(\mcl{D}^m)=m$. If $\Theta$ is a  $\mathcal{D}$-submodule of $\mathcal{R}^m$ and $\Delta$ is a submodule of $\Theta$, we say that 
$\Delta$ \textit{is primitive in} $\Theta$ if any submodule of $\Theta$ containing $\Delta$ and having rank equal to
$\rank(\Delta)$  is equal to $\Delta$. We see that the 
set of all nonzero primitive submodules of a fixed $\mathcal{D}$-submodule $\Theta$ of $\mathcal{R}^m$
is a partially ordered set with respect to set inclusion and its length is equal to $\rank(\Theta)$. When
$\Theta = \mathcal{D}^m$, we can even characterize the primitive submodules of $\mathcal{D}^m$ from the following observation:
\[\Delta\,\,\text{ is primitive }\,\,\Longleftrightarrow \Delta= K\Delta \cap \mathcal{D}^m \Longleftrightarrow \Delta= \mathcal{R}\Delta \cap \mathcal{D}^m\,.\]
This also shows that for any submodule $\Delta'$ of $\mathcal{D}^m$ there exists a unique primitive submodule $\Delta\supseteq \Delta'$ such that 
$\rank(\Delta)= \rank(\Delta')$, namely $\Delta := K\Delta'\cap \mathcal{D}^m$. \\

Let $\mathcal{R}$ have a topological ring structure in addition. We consider the topological group
$\GL(m,\mathcal{R})$ of $m\times m$ invertible matrices with entires in $\mathcal{R}$. It is obvious that any $g\in \GL(m,\mathcal{R})$ maps 
$\mathcal{D}$-submodules of $\mathcal{R}^m$ to $\mathcal{D}$-submodules of $\mathcal{R}^m$ preserving their rank and
inclusion relation. Let 
\[\frk{M}(\mcl{R},\mcl{D},m):= \{g\Delta\,:\, g\in \GL(m,\mcl{R}),\,\,\Delta \,\,\mbox{is a submodule of of}\,\,\mcl{D}^m\}\,.\]
We also denote the set of all nonzero primitive submodules of $\mcl{D}^m$, which is a poset of length $m$ with respect to inclusion relation 
as we have already seen, by $\frk{P}(\mcl{D},m)$.  \\

For a given function $||\cdot||: \frk{M}(\mcl{R},\mcl{D},m) \longrightarrow \mathbb{R}_{\geq 0}$, one says that $||\cdot||$ is 
\textit{norm-like} if the following three conditions hold: 
\begin{enumerate}[(N1)]
 \item \label{item:N1} For any $\Delta,\Delta'\in \frk{M}(\mcl{R},\mcl{D},m)$ with $\Delta'\subseteq \Delta$ and
 $\rank(\Delta')= \rank(\Delta)$, we always have $||\Delta'||\geq ||\Delta||$;
 \item \label{item:N2} there exists $C_{||\cdot||} \textgreater 0$ such that $||\Delta + \mcl{D}\gamma||\leq C_{||\cdot||} ||\Delta||\,||\mcl{D}\gamma||$
 holds for any $\Delta \in \frk{M}(\mcl{R},\mcl{D},m)$ and any $\gamma\notin \mcl{R}\Delta$; and 
 \item \label{item:N3} the function $\GL(m,\mcl{R})\longrightarrow \mathbb{R}_{\geq 0},\,g\mapsto ||g\Delta||$ is continuous for every submodule $\Delta$
 of $\mcl{D}^m$. 
\end{enumerate}

With the notations and terminologies defined so above, it is now time to state the `Quantitative nondivergence theorem'. 
\begin{theorem}\label{thm:qn}
 Let $B\subseteq X$ be a ball in a Besicovitch metric space $X$ and    $h:\tilde{B} \longrightarrow \GL(m,\mcl{R})$, 
 where $\tilde{B}:=3^m B$,  
be a continuous map. Suppose $\nu$ is a radon measure on $X$ which is D-Federer on $\tilde B$.  
Assume that a norm-like function $||\cdot||$ is given on $\frk{M}(\mcl{R},\mcl{D},m)$. Assume further that 
for some $C,\alpha \,\textgreater \,0$ and $\rho \in \displaystyle (0, 1/C_{||\cdot||}]$, the following conditions hold:
\begin{enumerate}[(C1)]
 \item \label{item:C1} for every $\Delta \in \frk{P}(\mcl{D},m)$, the function $x\mapsto ||h(x)\Delta||$ is $(C,\alpha)-good$ on $\tilde{B}$ w.r.t  $\nu$;
 \item \label{item:C2}for every $\Delta \in \frk{P}(\mcl{D},m)$, $\displaystyle \sup_{x\,\,\in B\cap \,\mbox{supp}\,(\nu)}||h(x)\Delta||\geq \rho$; and 
 \item \label{item:C3}$\forall x\in \tilde{B} \cap \mbox{supp}\,(\nu)$, 
 $\#\{\Delta \in \frk{P}(\mcl{D},m)\,:\,||h(x)\Delta||\,\textless \,\rho\} \,\textless \,\infty$\,. 
\end{enumerate}
Then for any positive $\varepsilon \leq \rho$, one has
\begin{equation}\label{eqn:qn}
\displaystyle  \nu\left(\left\{x\in B\,:\,||h(x)\gamma||\,\textless\, \varepsilon\,\,\mbox{for some}
\,\,\gamma \in \mcl{D}^m\setminus \{\mathbf{0}\}\right\}\right)\leq mC(N_X D^2)^m \left(\frac{\varepsilon}{\rho}\right)^{\alpha} \nu(B)\,,
\end{equation}
where $N_X$ is the `Besicovitch constant'.  

\end{theorem}
For the proof, 
see (\cite{KT}, \S 6, Theorem).

\subsection{The proof of Proposition \ref{main prop}}
From the definition of $K_{\varepsilon}$, as in $(\ref{eqn:cpt})$, it is obvious that for
$\mathbf{t}\in \mathfrak{a}^+$ and $\mathbf{x}\in B$, 
\[g_{\mathbf{t}}\overline{\tau}(\mathbf{f}(\mathbf{x}))\notin K_{\varepsilon} \Longleftrightarrow 
 ||(g_{\mathbf{t}}\tau(\mathbf{f}(\mathbf{x})))\mathbf{v}||\textless \varepsilon \,\,\mbox{for some}\,\,\mathbf{v}\in \Lambda^{n+1}
 \setminus \{\mathbf{0}\}\,.
\]
\noindent This inspires us to use  
Theorem \ref{thm:qn} in the setting
\[\begin{array}{rcl}\mcl{D}=\Lambda, \mcl{R}=F,
  X=F^d,m=n+1;\\ \nu, B, C, \alpha \,\,\mbox{and}\,\,D\,\,\mbox{as in Proposition \ref{main prop}};\\
  h(\mathbf{x})= g_{\mathbf{t}}\tau(\mathbf{f}(\mathbf{x}))\,\,\forall\mathbf{x}\in \tilde{B}\,;\,\end{array}
 \] and   $||\cdot||$ as the following:\\
 
 Since $\Lambda$ is a PID, any  submodule of the $\Lambda$ module $\Lambda^{n+1}$, being submodule of a free module of rank $n+1$, 
 is free of rank $\leq n+1$. Thus any nonzero $\Delta \in \frk{M}(F,\Lambda,n+1)$ has a $\Lambda$ basis, say 
 $\{\mathbf{v}_1,\,.\,. \,.\,,\mathbf{v}_j\}$, where $1\leq j\leq n+1$. We consider the $j$-vector 
 $\mathbf{w}:=\mathbf{v}_1\wedge \cdot\cdot\cdot \wedge \mathbf{v}_j\in \bigwedge^j (F^{n+1})$. Recall that 
 the $j$-vectors $e_{i_1}\wedge e_{i_2}\wedge \cdot\cdot\cdot\wedge e_{i_j}$
 with integers $1\leq i_1\textless i_2 \textless \cdot\cdot\cdot\textless i_j\leq n+1$ form a basis of $\bigwedge^j (F^{n+1})$ and
 thus $\bigwedge^j (F^{n+1})$ can be identified with $F^{\binom{n+1}{j}}$. Therefore one can naturally talk about 
 the \textit{supremum norm} on $\bigwedge^j (F^{n+1})$ using this identification. We  define 
 \[ ||\Delta||:= \text{ supremum norm of }\mathbf{w}\,.\]
  It is a routine 
 verification that this definition does not depend on the choice of the ordered basis 
 of $\Delta$. If $\Delta=\{\mathbf{0}\}$, we define $||\Delta||=1$. \\
 
 In order to prove that the just defined $||\cdot||$ is indeed norm-like, we need to verify the 
 conditions (N\ref{item:N1})-(N\ref{item:N3}). (N\ref{item:N1}) and 
 (N\ref{item:N3}) follow easily from the basic properties 
 of exterior product, while (N\ref{item:N2}) can be proved by a verbatim repetition of the proof of Lemma 5.1 of \cite{KM} as follows. \\

 We claim that $C_{||\cdot||}$ can be taken as $1$. If $\Delta=\{\bf{0}\}$ then it is immediate. Otherwise let
 $\{\mathbf{v}_1,\,.\,. \,.\,,\mathbf{v}_j\}$ be a basis of $\Delta$. Clearly $\{\mathbf{v}_1,\,.\,. \,.\,,\mathbf{v}_j,\gamma\}$ is a basis 
 of $\Delta+\Lambda \gamma$. Now writing 
 $\mathbf{v}_1\wedge \cdot\cdot\cdot \wedge \mathbf{v}_j= \displaystyle\sum_{\tiny \begin{array}{rcl}I\subseteq \{1,2,...,n+1\},\\ \# I=j\end{array}}w_I e_I$ and 
 $\gamma=\displaystyle \sum_{i=1}^{n+1} w_ie_i$ (in usual notations) and using the ultrametric property, we see that 
 \[\begin{array}{rcl}
    ||\Delta+\Lambda \gamma||=\left|\left|\displaystyle\sum_{{\tiny \begin{array}{rcl}I\subseteq \{1,2,...,n+1\},\\ \# I=j\end{array}}}w_I e_I\wedge \sum_{i=1}^{n+1} w_ie_i\,\right|\right|
    \leq
     \displaystyle \max_{1\leq i\leq n+1} \left|\left|\displaystyle\sum_{{\tiny \begin{array}{rcl}I\subseteq \{1,2,...,n+1\},\\ \# I=j\end{array}}}w_I w_i (e_I\wedge  e_i)\,\right|\right|

    \\ \leq \displaystyle \max_{1\leq i\leq n+1} \max_{{\tiny \begin{array}{rcl}I\subseteq \{1,2,...,n+1\},\\\# I=j\end{array}}}|w_Iw_i|\\ \leq 
     \displaystyle \max_{{\tiny \begin{array}{rcl}I\subseteq \{1,2,...,n+1\},\\ \#I=j\end{array}}}|w_I| \displaystyle \max_{1\leq i\leq n+1}|w_i|\\=
     ||\Delta||\,||\Lambda \gamma||\,.
    \end{array}
\]

Now we have to check the conditions (C\ref{item:C1}), (C\ref{item:C2}) and (C\ref{item:C3}) of  Theorem \ref{thm:qn}. From the discreteness 
of $\bigwedge^j (\Lambda^{n+1})$ in $\bigwedge^j (F^{n+1})$ for all $j=1,2,...,n+1$, (C\ref{item:C3}) is immediate. To 
investigate the validity of others, we have 
to do the explicit computation exactly in the similar manner to that of \S 3.3 in \cite{KW1}. \\

$\bullet$ \textbf{Checking (C\ref{item:C1}):} Here, for the sake of convenience in computation, 
it is customary to bring a few minor changes in some of the notations we have been using so far. For the rest of this section,
we write $\{\mathbf{e}_0,\mathbf{e}_1,...,\mathbf{e}_n\}$  the standard basis of $F^{n+1}$  
 and for
\begin{equation} \label{eqn:I} I=\{i_1,...,i_j\}\subseteq \{0,...,n\}\text{ where \,}i_i\,\textless \,i_2\,\textless \cdot\cdot\cdot\textless \,i_j\,,\end{equation}
we let $\mathbf{e}_I$ denote $\mathbf{e}_{i_1}\wedge \cdot \cdot \cdot \wedge \mathbf{e}_{i_j}$. Similarly, it will be convenient to put
any $\mathbf{t}\in \frk{a}^+$ as 
\[\mathbf{t}=(t_0,t_1,...,t_n)\,\,\,\mbox{where}\,\,t_0=\displaystyle \sum_{i=1}^n t_i\,.\]
Let us observe that for any $\mathbf{y}\in F^n$, $\tau(\mathbf{y})$ fixes $\mathbf{e}_0$  and sends any other $\mathbf{e}_i$ 
to $\mathbf{e}_i+y_i\mathbf{e}_0$. Thus for any $I$ as in (\ref{eqn:I}), we have 
\begin{equation}\label{eqn:5.1}
 \tau(\mathbf{y})\mathbf{e}_I= \left \{\begin{array}{rcl}
                                        \mathbf{e}_I & \mbox{if}\,\,0\in I\\
                                        \mathbf{e}_I + \sum_{i\in I} \pm y_i \mathbf{e}_{I\cup \{0\}\setminus \{i\}} & \mbox{otherwise}\,.
                                       \end{array}
 \right.  
\end{equation}
Likewise, we can also see that for any $I$ as in (\ref{eqn:I}), 
\begin{equation}\label{eqn:5.2}
 g_{\mathbf{t}}\mathbf{e}_I= \left \{\begin{array}{rcl}
                                        \displaystyle T^{t_0- \sum_{i\in I\setminus \{0\}}t_i} \,\mathbf{e}_I & \mbox{if}\,\,0\in I\\
                                        T^{-\sum_{i\in I} t_i}\,\mathbf{e}_I & \mbox{otherwise}\,.
                                       \end{array}
 \right. 
\end{equation}

\indent Suppose $\Delta \in \frk{P}(\Lambda,n+1)$  and $\{\mathbf{v}_1,\,.\,. \,.\,,\mathbf{v}_j\}$ is a basis of $\Delta$
and let \[\mathbf{w}:= \mathbf{v}_1\wedge \cdot\cdot\cdot \wedge \mathbf{v}_j= \displaystyle\sum_{\tiny \begin{array}{rcl}I
\subseteq \{0,...,n\},\\ \# I=j\end{array}}w_I \mathbf{e}_I\,;\,\,\,\,w_I\in \Lambda\,.\] From (\ref{eqn:5.1}) and (\ref{eqn:5.2}), it 
follows that for any $\mathbf{x}\in \tilde B$, one has 
\[h(\mathbf{x})\mathbf{w}= \displaystyle\sum_{\tiny \begin{array}{rcl}I\subseteq \{0,...,n\},\\ \# I=j\end{array}} h_I(\mathbf{x})\mathbf{e}_I\,\] 
where 
\begin{equation} \label{eqn:coeff}  h_I(\mathbf{x}):= \left \{ \begin{array}{rcl}\displaystyle T^{- \sum_{i\in I} t_i}w_I & \mbox{if} \,\, 0\notin I \\
                         \displaystyle T^{\sum_{i\notin I} t_i}(w_I + \sum_{i\notin I} \pm w_{I \cup \{i\} \setminus \{0\}} f_i(\mathbf{x}))           &\mbox{otherwise} .
                      
                        \end{array} \right.\end{equation}
 
In particular, the coordinate maps $h_I$ of the map $\mathbf{x}\mapsto  h(\mathbf{x})\mathbf{w},\,\mathbf{x}\in \tilde B$ are $F$-linear combinations 
of $1,f_1,\,.\,.\,.\,,f_n$ and hence, by (i) of Proposition \ref{main prop}, all of them 
are $(C,\alpha)-good$ on $\tilde B$ with respect to $\nu$. Therefore, from (\ref{item:sup}) of Lemma \ref{lem:C,alpha}, it follows that the 
function $$\displaystyle \mathbf{x}\mapsto ||h(\mathbf{x})\Delta||=  ||h(\mathbf{x})\mathbf{w}||=
\max_{I}|h_I(\mathbf{x})|$$ is $(C,\alpha)-good$ on $\tilde{B}$ with respect to  $\nu$. 
Thus (C\ref{item:C1}) is established.\\

$\bullet$ \textbf{Checking (C\ref{item:C2}):} Let $\Delta \in \frk{P}(\Lambda,n+1)$, $\{\mathbf{v}_1,\,.\,. \,.\,,\mathbf{v}_j\}$ be
a basis of $\Delta$
and let \[\mathbf{w}:= \mathbf{v}_1\wedge \cdot\cdot\cdot \wedge \mathbf{v}_j= \displaystyle\sum_{\tiny \begin{array}{rcl}I
\subseteq \{0,...,n\},\\ \# I=j\end{array}}w_I \mathbf{e}_I\,;\,\,\,\,w_I\in \Lambda\,.\] 
\textit{Case 1:} Assume $w_I=0$ whenever $0\notin I$. Then there must be  some $J\subseteq \{0,...,n\}$  containing $0$ such that 
$w_J\neq 0$ as all $w_I$ 
can not be zero.  Pick any $\mathbf{t}\in \frk{a}^+$.  Now 
from (\ref{eqn:coeff}),  we see that 
$$|h_J(\mathbf{x})|=\displaystyle |T^{\sum_{i\notin J} t_i}w_J|\geq1\text{ \,\,for any }\mathbf{x}\in \tilde B\,.$$  Therefore in this case, we have
\begin{equation}\label{eqn:C1}\begin{array}{rcl} \displaystyle \sup_{\tiny \mathbf{x}\in B\,\cap \,\mbox{supp}\,(\nu)}||h(\mathbf{x})\Delta||
=\sup_{\tiny \mathbf{x}\in B\,\cap \,\mbox{supp}\,(\nu)}||h(\mathbf{x})\mathbf{w}|| = \displaystyle 
\sup_{\tiny \mathbf{x}\in B\,\cap \,\mbox{supp}\,(\nu)} \max_{I} |h_I(\mathbf{x})| \\ \displaystyle \geq \sup_{\tiny \mathbf{x}\,\,\in B\,\cap \,\mbox{supp}\,(\nu)} |h_J(\mathbf{x})|\geq 1
\\ \forall \,\mathbf{t}\in \frk{a}^+\,.\end{array}\end{equation}
\textit{Case 2:} Suppose $w_I\neq 0$ for some $I\subseteq \{1,...,n\}$. Choose $l\in\{1,...,n\}$  such that $t_l= \max_{1\leq i\leq n}t_i$.
If $l\in I$, set $J= I\cup \{0\}\setminus \{l\}$.  Clearly
$J$ contains $0$ but does not contain $l$. In view of (\ref{eqn:coeff}),
the  coefficient of $f_l$  
in the expression  of $h_J$ is easily seen to be $\pm T^{\sum_{i\notin J} t_i}w_I$ and its absolute value is

\begin{equation}\label{eqn:lbc}\displaystyle |T^{\sum_{i\notin J} t_i}w_I|
\geq e^{\sum_{i\notin J} t_i}\geq e^{t_l}\geq e^{t_0 /n}=e^{||\mathbf{t}||/n}\,.\end{equation}
If $l \notin I$,  choose any $i\in I$ and let $J= I\cup \{0\}\setminus \{i\}$. Like before, 
$\pm T^{\sum_{i\notin J} t_i}w_I$ turns out as the coefficient of $f_i$ in $h_J$ so that we obviously  get the analogue of (\ref{eqn:lbc}). Thus in this case,
there always exists $J$ such that \\
\begin{equation}\label{eqn:coeff*} \mbox{at least one 
of the coefficients of}\,\,\, f_1,f_2,...\,,f_n \,\,\mbox{in}\,\, h_J\,\,\mbox{has absolute value}\geq e^{||\mathbf{t}||/n}\,.\end{equation}

Now, from the assumption (ii) of Proposition \ref{main prop}, it follows that there exists $\delta\, \textgreater\, 0$ such that 
$\sup_{\tiny \mathbf{x}\in B\,\cap\, \mbox{supp}\,(\nu)} |c_0+c_1f_1(\mathbf{x})+\cdot \cdot \cdot +c_nf_n(\mathbf{x})|\geq \delta$ for any $c_0,c_1,...,c_n\in F$ with
$\max_{\tiny 0\leq i\leq n} |c_i|\geq 1$. 
We choose $M\in \mathbb{N}$ such that $\delta e^M\geq 1$.  \\

Let $||\mathbf{t}||\geq nM$. Then, because of (\ref{eqn:coeff*}), one surely has
at least one 
of the coefficients of $f_1,f_2,...\,,f_n$ in $\frac{1}{T^M}\,h_J$ has absolute value at least $1$ and thus 

\[\sup_{\tiny \mathbf{x}\,\cap \,\mbox{supp}\,(\nu)} \left|\frac{1}{T^M}\,h_J (\mathbf{x})\right|\geq \delta\,. \] This gives, 
\[\sup_{\tiny \mathbf{x}\,\cap \,\mbox{supp}\,(\nu)} \left|h_J (\mathbf{x})\right|\geq \delta e^M\geq 1\,.\] 
So, even here, we can see that
\begin{equation}\label{eqn:C2}
 \begin{array}{rcl}\displaystyle \sup_{\tiny \mathbf{x}\,\,\in B\,\cap \,\mbox{supp}\,(\nu)}||h(\mathbf{x})\Delta||
 = \sup_{\tiny \mathbf{x}\,\,\in B\,\cap \,\mbox{supp}\,(\nu)}||h(\mathbf{x})\mathbf{w}||=
 \sup_{\tiny \mathbf{x}\,\,\in B\,\cap \,\mbox{supp}\,(\nu)} \max_{I} 
|h_I(\mathbf{x})|\\ \displaystyle \geq \sup_{\tiny \mathbf{x}\,\,\in B\,\cap \,\mbox{supp}\,(\nu)} |h_J(\mathbf{x})|\geq 1 \\ \forall\, ||\mathbf{t}||\geq nM\,.
\end{array} \end{equation} 

\noindent Letting $\rho=1$, $(C2)$ is thus immediate from (\ref{eqn:C1}) and (\ref{eqn:C2})  whenever $||\mathbf{t}||\geq nM$.\\

Finally, $\tilde C$ and $s$ are taken as $(n+1)CD^{2(n+1)}$ and $nM$ repectively, and 
one  applies Theorem \ref{thm:qn} to show (\ref{eqn:prop}). \,\,$\Box$. 

\section{Explicit constants: an example}\label{section:example}
In this section, we  talk about a simple application 
of our Theorem (\ref{thm:1.2}) to a concrete example, with special attention on the explicit constant $\varepsilon_0$. For us, here $d=1, n=2$ and $\nu$ is the unique Haar measure on $F$ that satisfies 
$\nu(B[0;1])=1$. It is not difficult to show that $\nu$ is $e^2$-Federer. Let 
\[f: B(0,1)\longrightarrow F^2\,,x\mapsto (x,x^2)\,.
\] We claim that $f$ is $(2,1/2)-good$, i.e. in other words, so is any $\phi\in F[x]$ 
having degree $\leq 2$. To see this, we shall 
apply the same technique
which used in the proof of proposition 3.2 of \cite{KM}.\\

Let $\varepsilon\,\textgreater \,0$ and $\mathfrak{B}\subseteq B(0;1)$. We have to show 
\begin{equation}\label{eqn:2,1/2}
\nu(\{x\in \mathfrak{B}\,:\, |\phi(x)|\,\textless\, \varepsilon\})\leq 2 \left(\frac{\varepsilon}{||\phi||_{\mathfrak{B}}}\right)^{1/2} \nu(\mathfrak{B})\,.
\end{equation}
For convenience, put $S:= \{x\in \mathfrak{B}\,:\, |\phi(x)|\,\textless \,\varepsilon\}$. If $\nu(S)$, 
i.e. the LHS of (\ref{eqn:2,1/2}), is $0$ then there is nothing to prove. Otherwise, we will show that 
\[m\leq 2 \left(\frac{\varepsilon}{||\phi||_{\mathfrak{B}}}\right)^{1/2} \nu(\mathfrak{B})\] of equivalently,
\begin{equation}\label{eqn:2,1/2, suff}
 |\phi(x)|\leq \varepsilon \left(\frac{\nu(\frk{B})}{m/2}\right)^2 \text{\,\,\,for all } x\in \frk{B}\,,
\end{equation}  whenever  $0\,\textless \,m\,\textless \,\nu(S)$.\\

From the continuity of $\phi$, we see that for each $x\in S$, there is 
a ball $B_x$ with center at $x$ and radius $ \textless \, \frac{m}{2}$ such that $B_x\subseteq S$.
Now from the Besicovitch nature of $F$, one can extract a countable subcover $\{B_1,B_2,...\}$  consisting of mutually disjoint open balls 
from the cover $\{B_x\,:x\in S\}$ of $S$. Clearly $\nu(B_i)\leq \frac{m}{2}$ for each $i$. Thus in view of their size,
it follows that the subcover has at least three balls. Let us denote their centers as $x_1, x_2\text{ and } x_3$. Then the centers $x_i\in S$ and they must 
satisfy 
\begin{equation}\label{eqn:center}
 |x_i- x_j|\geq \frac{m}{2}\text{\,\,\, for all }i,j= 1,2,3\,;\, i\neq j\,.
\end{equation}
It is now time to employ the `Lagrange's interpolation formula' to complete the proof. By this formula, we can write $\phi$ as
\begin{equation}\label{eqn:lag}\begin{array}{rcl}\displaystyle \phi(x)=\phi(x_1)\frac{(x-x_2)(x-x_3)}{(x_1-x_2)(x_1-x_3)}+\phi(x_2)\frac{(x-x_1)(x-x_3)}{(x_2-x_1)(x_2-x_3)}\\ \displaystyle +\phi(x_3)
\frac{(x-x_1)(x-x_2)}{(x_3-x_1)(x_3-x_2)}\,.\end{array}
 \end{equation}
 As $\frk{B}$ is a ball, certainly there exist $a\text{ and }m\in \mathbb{N}$ such that $\frk{B}=B[a;\frac{1}{e^m}]$. Therefore 
 $\diameter (\frk{B})=\frac{1}{e^m}=\nu(\frk{B})$. In view of this, (\ref{eqn:center}) and (\ref{eqn:lag}), it follows at once that
 \[|\phi(x)|\leq \varepsilon\frac{(\diameter (\frk{B}))^2}{m^2/4}= \varepsilon \left(\frac{\nu(\frk{B})}{m/2}\right)^2 \text{\,\,\,for all } x\in \frk{B}\,;\]
and that shows (\ref{eqn:2,1/2, suff}). \\

Thus all the conditions of the hypothesis of Theorem \ref{thm:1.2} hold here and so that existence of desired $\varepsilon_0\, \textgreater\, 0$ is 
confirmed. We are interested to compute it. In the proof of Theorem \ref{thm:1.2}, we have also observed that 
our $\varepsilon_0$ can be taken as any positive quantity which is $\displaystyle \textless \,\frac{1}{\tilde C^{1/\alpha}}$; where 
$\tilde C$ was set, as we did in the proof of Proposition \ref{main prop}, as $(n+1)CD^{2(n+1)}$. Therefore in our example, 
we obtain that
\[\varepsilon_0 \textless\, \frac{1}{\tilde C^2}= \frac{1}{(3\times 2\times (e^2)^ {2\times 3})^2}=\frac{1}{36\,e^{24}}\,.\]

\section{Appendix: The density theorem}\label{section:density}  Fix $d\in \mathbb{N}$ and a radon measure $\nu$ on $F^d$. Let $\Omega$ be a measurable subset of $F^d$ and $\mathbf{x}\in \mbox{supp}\,(\nu)$. 
We say that $\mathbf{x}$ is a \textit{point of density} of $\Omega$ if 
$$ \displaystyle \lim_{\tiny \begin{array}{rcl}\nu (B)\rightarrow 0\\ \mathbf{x}\in B\end{array}}\frac{\nu(B\cap \Omega)}{\nu(B)}=1\,.$$ 
\noindent The above definition is nothing but the counterpart of the classical notion  `point of Lebsegue density' in our setting. Likewise, it is thus natural to expect the following:
\begin{theorem}\label{thm:3.1}
 Suppose $\Omega$ is a measurable subset of $F^d$. Then : 
 \begin{enumerate}
  \item Almost every $\mathbf{x}\in \Omega$ is a point of density of $\Omega$\,.
  \item Almost every $\mathbf{x}\notin \Omega$ is not a point of density of $\Omega$\,.
 \end{enumerate}

\end{theorem}

Now to prove Theorem \ref{thm:3.1}, we shall take up the same strategy of the proof for euclidean spaces. Namely, we 
develop  our version 
of the `Lebsegue differentiation theorem' first and Theorem \ref{thm:3.1} follows then as a consequence. To begin with, we need to introduce `locally integrable' functions as follows:\\

\indent A measurable function $f$ on $F^d$ is said to be \textit{locally integrable} if for every ball $B\subseteq F^d$, the 
function $f\chi_{B} \in L^1(F^d,\nu)$.\\

\noindent With this terminology, let us state our `differentiation theorem': 
\begin{theorem}\label{thm:3.2} Let $f$ be a locally integrable function on $F^d$. Then 
\begin{equation}\label{eqn:3.1}\displaystyle \lim_{\tiny \begin{array}{rcl}\nu (B)\rightarrow 0\\ \mathbf{x}\in B\end{array}} 
\frac{1}{\nu(B)} \int_{B} f\,\,d\nu= f(\mathbf{x})\,\,\mbox{for almost all}\,\,\mathbf{x}\in \mbox{supp}\,(\nu)\,.\end{equation}
                                                 
\end{theorem}An application of the above theorem to the characteristic function of $\Omega$ immediately yields Theorem \ref{thm:3.1}.  \\

To prove  Theorem \ref{thm:3.2}, let us observe that $f$ can be assumed to be integrable without any loss in generality. To see this, 
suppose that the conclusion of the theorem is established for integrable functions. Now for any $n \in \mathbb{N}$, we
apply (\ref{eqn:3.1}) to the integrable 
function $f\chi_{B(\mathbf{0};n)}$ and obtain that the set of all $\mathbf{x}\in \mbox{supp}\,(\nu)\cap B(\mathbf{0};n)$ for which 
\begin{equation}\label{eqn:3.2}\displaystyle \lim_{\tiny \begin{array}{rcl}\nu (B)\rightarrow 0\\ \mathbf{x}\in B\end{array}} 
\frac{1}{\nu(B)} \int_{B} f\,\,d\nu \neq f(\mathbf{x})\end{equation}
has zero measure. Since any $\mathbf{x}\in \mbox{supp}\,(\nu)$ that satisfies (\ref{eqn:3.2}) is contained in some $B(\mathbf{0};n)$, we are done. Thus 
we shall always let $f\in L^1(F^d,\nu)$ in the rest of this section.\\

It suffices to show that, for all $\alpha\, \textgreater \,0$, 
$$E_{\alpha} := \left \{\mathbf{x}\in \mbox{supp}\,(\nu)\,:\,\displaystyle \limsup_{\tiny \begin{array}{rcl}\nu (B)\rightarrow 0\\ 
\mathbf{x}\in B\end{array}}\left|\frac{1}{\nu(B)}\int_{B} f\,\,d\nu\,-f(\mathbf{x})\right|\,\textgreater \,2\alpha \right\}$$ is null 
with respect to $\nu$.
To achieve this, we need to introduce the `Hardy-Littlewood maximal function' and make use of its main 
property. The relevant definition goes as follows.\\

Let $f\in  L^1(F^d,\nu)$. The \textit{maximal function} of $f$, denoted by $f^*$, is defined as 
$$f^{*}(\mathbf{x})= \displaystyle \sup_{\mathbf{x}\in B}\, \frac{1}{\nu(B)}\int_B |f|\,\,d\nu\,\,\,for\,\,\mathbf{x}\in\mbox{supp}\,(\nu)\,\,
\mbox{and}\,\, \infty \,\,\,\mbox{otherwise} \,.$$
It is easy to see that, $\forall \alpha \in \mathbb{R}$, the set $\{\mathbf{x}\in F^d\,:\,f^{*}(\mathbf{x})\,\textgreater \,\alpha \}$ is
open,
because if $\mathbf{x}\in \mbox{supp}\,(\nu)$ and $f^{*}(\mathbf{x})\,\textgreater \,\alpha$ then there exists a ball $B$ containing $\mathbf{x}$ for which 
$$\frac{1}{\nu(B)}\int_B |f|\,\,d\nu \,\textgreater\, \alpha\,.$$ Clearly for any $\overline{\mathbf{x}} \in B \cap \mbox{supp}\,(\nu)$, one has
$$f^*(\overline{\mathbf{x}})\geq \frac{1}{\nu(B)}\int_B |f|\,\,d\nu\, \textgreater \,\alpha\,.$$ Hence $f^*$ is measurable. The main
property of this maximal function is given by the following theorem.
\begin{theorem}\label{thm:3.3}
 Let $f$ be integrable. Then for all $\alpha \,\textgreater \,0$
 \begin{equation} \label{eqn:3.3}
 \nu (\{\mathbf{x}\in F^d\,:\, f^*(\mathbf{x})\,\textgreater\, \alpha\}) \leq \frac{1}{\alpha}\,||f||_1\,, 
 \end{equation}
where $||f||_1=\int_{F^d} |f|\,\,d\nu$. As a consequence, $f^*(\mathbf{x})\,\textless \,\infty$ for almost all $\mathbf{x}$.
\end{theorem}
\textit{Proof} : Since $\nu$ is radon and the 
set $\mathcal{A}_{\alpha} := \{\mathbf{x}\in F^d\,:\, f^*(\mathbf{x})\,\textgreater \,\alpha\}$ is open, as seen earlier, one has 
$$\nu(\mathcal{A}_{\alpha})=\displaystyle \sup_{\tiny \begin{array}{rcl}K\subseteq \mathcal{A}_{\alpha}\\K\,\mbox{ compact}\end{array}} \nu(K)\,.$$
 It is thus enough to show that, for any compact subset $K$ of $\mathcal{A}_{\alpha}$,
$$\nu(K)\,\textless\, \frac{1}{\alpha}\,||f||_1\,.$$ For each $\mathbf{x}\in K \cap \mbox{supp}\,(\nu)$, we have a ball 
$B_{\mathbf{x}}$ satisfying 
$$\frac{1}{\nu(B_{\mathbf{x}})}\int_{B_{\mathbf{x}}} |f|\,\,d\nu\, \textgreater \,\alpha\,,$$ i.e. 
\begin{equation}\label{eqn:3.4}\frac{1}{\alpha}\int_{B_{\mathbf{x}}} |f|\,\,d\nu \,\textgreater \,\nu(B_{\mathbf{x}})\,.\end{equation} The compact set 
$ K \cap \,\mbox{supp}\,(\nu)$ can be covered by finitely many such balls, say $B_1\,\,B_2\,,\,.\,.\,.\,,\,B_r$. Without any loss in generality, 
we can assume that the collection $\{B_i\}_{i=1}^r$ of balls is mutually disjoint.  \\

Now, in view of (\ref{eqn:3.4}), we find that 
$$ \begin{array}{rcl}\nu (K)= \nu (K \cap \,\mbox{supp}\,(\nu))\leq \displaystyle \sum_{i=1}^{r} \nu(B_{i})\, 
\textless \,\displaystyle \sum_{i=1}^{r} \frac{1}{\alpha}\int_{B_{i}} |f|\,\,d\nu = 
\displaystyle \frac{1}{\alpha} \sum_{i=1}^{r} \int_{B_{i}} |f|\,\,d\nu \\ = \displaystyle \frac{1}{\alpha} \int_{\cup_{i=1}^r B_i} |f|\,\,d\nu
\\ \leq \displaystyle \frac{1}{\alpha}\int_{F^d} |f|\,\,d\nu\\= \displaystyle \frac{1}{\alpha} ||f||_1\,.\,\,\,\,\,\,\Box \end{array}$$
\indent We shall also need the following lemma 
\begin{lemma}\label{lem:3.1}
 Let $g$ be an integrable function on $F^d$ and $\mathbf{x}\in \mbox{supp}\,(\nu)$ is a point of continuity of $g$.  Then 
$$\displaystyle \lim_{\tiny \begin{array}{rcl}\nu (B)\rightarrow 0\\ \mathbf{x}\in B\end{array}} 
\frac{1}{\nu(B)} \int_{B} g\,\,d\nu= g(\mathbf{x})\,.$$
\end{lemma}
\textit{Proof} : Let $\varepsilon\, \textgreater \,0$. From the continuity of $g$ at the
point $\mathbf{x}$, we get $r\,\textgreater \,0$ such that
$|g(\mathbf{y})-g(\mathbf{x})|\,\textless\, \varepsilon$ for each $\mathbf{y}\in B(\mathbf{x};r)$. 
We set $$\delta := \frac{1}{2} \,\nu(B(\mathbf{x};r))\,.$$ Pick any ball 
$B$ containing $\mathbf{x}$ with $\nu(B)\,\textless\, \delta$. Clearly $B\subseteq B(\mathbf{x};r)$, because 
otherwise $B(\mathbf{x};r) \subsetneqq B$ will hold, as $\mathbf{x}\in B\cap B(\mathbf{x};r)$. But then we will have
$\delta\, \textless\, \nu (B(\mathbf{x};r))\leq \nu(B)$ which is impossible. Now, it is easy to see  that
$$\left|\frac{1}{\nu(B)} \int_{B} g\,\,d\nu- g(\mathbf{x})\right|\leq \frac{1}{\nu(B)} \left|\int_{B} (g-g(\mathbf{x}))\,d\nu \right|\leq\frac{1}{\nu(B)} \int_{B} \left|g-g(\mathbf{x})\right|\,d\nu\leq 
\frac{\varepsilon \nu(B)}{\nu(B)}= \varepsilon\,.\,\,\,\,\Box$$\\

Let us prove Theorem \ref{thm:3.2} now. Assume $\varepsilon\,\textgreater\, 0$. Since $F^d$ is locally compact and hausdorff, 
$C_c(F^d)$ is dense in $L^1(F^d,\nu)$. So we can choose a continuous and compactly supported function $g$ such that 
\begin{equation}\label{eqn:3.5}
 ||f-g||_1 \,\textless \,\varepsilon\,.
\end{equation}

For any $\mathbf{x}\in \mbox{supp}\,(\nu)$, writing $\frac{1}{\nu(B)}\int_{B} f\,\,d\nu\,-f(\mathbf{x})$ as 
$$\frac{1}{\nu(B)}\int_{B} (f-g)\,\,d\nu+ \frac{1}{\nu(B)} \int_{B} g\,\,d\nu\,- g(\mathbf{x})+(g(\mathbf{x})-f(\mathbf{x}))\,,$$
we can see that 
\begin{equation}\label{eqn:3.6}
 \begin{array}{rcl}\displaystyle \limsup_{\tiny \begin{array}{rcl}\nu (B)\rightarrow 0\\ 
\mathbf{x}\in B\end{array}}\left|\frac{1}{\nu(B)}\int_{B} f\,\,d\nu\,-f(\mathbf{x})\right|\leq \displaystyle \limsup_{\tiny \begin{array}{rcl}\nu (B)\rightarrow 0\\ 
\mathbf{x}\in B\end{array}}\frac{1}{\nu(B)}\int_{B} |f-g|\,\,d\nu \\ + \displaystyle \limsup_{\tiny \begin{array}{rcl}\nu (B)\rightarrow 0\\ 
\mathbf{x}\in B\end{array}}\frac{1}{\nu(B)}\left| \int_{B} g\,\,d\nu\,- g(\mathbf{x})\right|\\ + |g(\mathbf{x})-f(\mathbf{x})|\,.\end{array}
\end{equation}
Let us observe that first term on the RHS of (\ref{eqn:3.6}) is $\leq (f-g)^*(\mathbf{x})$ and the middle 
term vanishes in view of Lemma \ref{lem:3.1}. Thus we have

\[
 \displaystyle \limsup_{\tiny \begin{array}{rcl}\nu (B)\rightarrow 0\\ 
\mathbf{x}\in B\end{array}}\left|\frac{1}{\nu(B)}\int_{B} f\,\,d\nu\,-f(\mathbf{x})\right|\leq (f-g)^*(\mathbf{x}) + 
|g(\mathbf{x})-f(\mathbf{x})|\,.
\]

\noindent This shows that $E_{\alpha}\subseteq F_{\alpha}\,\cup \,G_{\alpha}$, where 
$$F_{\alpha}:= \{\mathbf{x}\in F^d\,:\, (f-g)^*(\mathbf{x})\,\textgreater\, \alpha\}\,\,\mbox{and}\,\,G_{\alpha}:=\{\mathbf{x}\in F^d\,:\,
|g(\mathbf{x})-f(\mathbf{x})|\,\textgreater \,\alpha\}\,.$$
We shall now estimate $\nu(F_{\alpha})$ and $\nu(G_{\alpha})$ one by one. On the one hand, as $f-g$ is integrable, it is immediate  that
\begin{equation}\label{eqn:3.7}
 \nu(G_{\alpha})\leq \frac{1}{\alpha}\,||f-g||_1\,. 
\end{equation}
On the other hand, Theorem \ref{thm:3.3} provides 
\begin{equation}\label{eqn:3.8}
 \nu(F_{\alpha})\leq \frac{1}{\alpha}\,||f-g||_1\,.
\end{equation}
Finally from (\ref{eqn:3.7}), (\ref{eqn:3.8}) and (\ref{eqn:3.5}), it follows that
$$\nu(E_{\alpha})\leq \nu(F_{\alpha})+\nu(G_{\alpha})\leq 2\,\frac{1}{\alpha}\,
||f-g||_1\leq \frac{2}{\alpha}\,\varepsilon\,.\,\,\,\,\Box$$


\begin{thebibliography}{99}

\bibliographystyle{aplha}

\bibitem{AGP} J. S. Athreya, A. Ghosh and A.Prasad, \textit{Ultrametric logarithm laws, II}, Monatsh. Math. 167 (2012), no. 3, 333--356.

\bibitem{AB}  B. Adamczewski,  Y. Bugeaud,  \textit{On the Littlewood conjecture in fields of power series}, Probability and number theory-Kanazawa 2005, 1--20, Adv. Stud. Pure Math., 49, Math. Soc. Japan, Tokyo, 2007.

\bibitem{Ba1} R.C. Baker, \textit{Metric diophantine approximation on manifolds}, J. Lond. Math. Soc. (2) 14 (1976), 43--48.

\bibitem{Ba2} \bysame, \textit{DirichletÕs theorem on diophantine approximation}, Math. Proc. Cambridge Phil. Soc. 83 (1978), 37--59.

\bibitem{Bu} Y. Bugeaud, \textit{Approximation by algebraic integers and Hausdorff dimension}, J. London Math. Soc. (2) 65 (2002), no. 3, 547--559.

\bibitem{DS1} H. Davenport and W.M. Schmidt, \textit{Dirichlet's theorem on diophantine approximation}, in: Symposia Mathematica, Vol. IV (INDAM, Rome, 1968/69), pp. 113--132, 1970.

\bibitem{DS2} \bysame, \textit{Dirichlet's theorem on diophantine approximation. II}, Acta Arith. 16 (1969/1970) 413--424.

\bibitem{deM} B. deMathan, \textit{Approximations diophantiennes dans un corps local}, Bull. Soc. Math. France 21 (Suppl. MŽm.) (1970)
1--93.

\bibitem{DRV} M. Dodson, B. Rynne, and J. Vickers, \textit{DirichletÕs theorem and Diophantine approximation on manifolds}, J. Number Theory 36 (1990), no. 1, 85--88.

\bibitem{DKJ} M. M. Dodson, S. Kristensen and J. Levesley, \textit{A quantitative Khintchine- Groshev type theorem over a field of formal series}, to appear in Indag. Math. (N.S.).

\bibitem{G-pos}
A. Ghosh, \textit{Metric Diophantine approximation over a local field of positive characteristic}, J. Number Theory 124 (2007), no. 2, 454--469.

\bibitem{GR} A. Ghosh and R. Royals, \textit{An extension of Khintchine's theorem}, Acta Arithmetica, 167 (2015), 1--17.

\bibitem{KM} 
D. Y. Kleinbock and G. A. Margulis, \textit{Flows on Homogeneous Spaces and Diophantine Approximation on Manifolds}, Annals of Mathematics, \textbf{148} (1998), 339--360.

\bibitem{KLW} 
D. Kleinbock, E. Lindenstrauss and B. Weiss, \textit{On fractal measures and Diophantine approximation}, Selecta Math. \textbf{10} (2004), 479--523.

\bibitem{KT} 
D. Kleinbock and G. Tomanov, \textit{Flows on $S$-arithmetic homogeneous spaces and 
applications to metric Diophantine approximation}, Comm. Math. Helv. \textbf{82} (2007), 519--581. 

\bibitem{KT-pre}
D. Kleinbock and G. Tomanov, MPI preprint.

\bibitem{KW1} 
D. Kleinbock and B. Weiss, \textit{Dirichlet's theorem on diophantine approximation and homogeneous flows}, J. Mod. Dyn. \textbf{4} (2008), 43--62. 

\bibitem{KW2}
\bysame, \textit{Friendly measures, homogeneous flows and singular vectors}, in: Algebraic and Topological Dynamics, Contemp. Math. 211, Amer. Math. Soc., Providence, RI, 2005, pp. 281--292.

\bibitem{Kris}
Simon Kristensen, \textit{On well approximable matrices over a field of formal series}, Math. Proc. Camb. Phil. Soc. (2003), 135(2), 255--268.

\bibitem{Las1}  A. Lasjaunias,  \textit{A survey of Diophantine approximation in fields of power series}, Monatsh. Math. 130 (2000), no. 3, 211--229.

\bibitem{Las2} A. Lasjaunias, \textit{Diophantine approximation and continued fractions in power series fields}, Analytic number theory, 297--305, Cambridge Univ. Press, Cambridge, 2009. 

\bibitem{Mahler}
Kurt Mahler, \textit{An analogue to MinkowskiÕs Geometry of numbers in a field of series}, Ann. Math., 2nd Ser. 1941, Vol.42, No.2. 488--522.

\bibitem{Mat}
P. Mattila, \textbf{Geometry of sets and measures in Euclidean space.
Fractals and rectifiability}, Cambridge Stud. Adv. Math. \textbf{44} Cambridge
University Press, Cambridge, 1995.

\bibitem{Schi}
W.H.Schikhof, \textit{Ultrametric Calculus, an introduction to p-adic analysis}, Cambridge studies in advanced mathematics, 4. Cambridge University Press, (1984).

\bibitem{Shah1} N. A. Shah, \textit{Equidistribution of expanding translates of curves and Dirichlet's theorem on Diophantine approximation},
Invent. Math. 177 (2009), no. 3, 509--532. 

\bibitem{Shah2} N. A. Shah, \textit{Expanding translates of curves and Dirichlet-Minkowski theorem on linear forms},  
J. Amer. Math. Soc. 23 (2010), no. 2, 563--589.

\bibitem{Shah3} N. A. Shah, \textit{Equidistribution of translates of curves on homogeneous spaces and Dirichlet's approximation}, 
Proceedings of the International Congress of Mathematicians. Volume III, 1332--1343, Hindustan Book Agency, New Delhi, 2010. 

\bibitem{Spr1} V.G. Sprindzuk, MahlerÕs Problem in Metric Number Theory, Transl. Math. Monogr., vol. 25, Amer. Math. Soc.,
1969 (translated from the Russian by B. Volkmann).

\bibitem{Spr2} V.G. Sprindzuk, Achievements and problems in Diophantine approximation theory, Russian Math. Surveys 35
(1980) 1--80.


\end{thebibliography}
\end{document}